\documentclass[a4paper]{article}
\usepackage{a4}
\usepackage{amssymb}
\usepackage{graphicx}
\begin{document}

\title{Flashes of noncommutativity}
\author{Alejandro Rivero\thanks{Zaragoza University at Teruel.  
           {\tt arivero@unizar.es}}}

\maketitle

\begin{abstract}
Noncommutativity lays hidden in the proofs of classical dynamics. Modern
frameworks can be used to bring it to light:  *-products, groupoids,
q-deformed calculus, etc. 
\end{abstract}

\section*{Flash one.}

Time ago, newborn Classical Mechanics simply described how
 the inertial law was 
disturbed
by the action of a force. One can consider two inertial trajectories from
$x$ to $y$ and then from $y$ to $z$,  where a change is applied in a 
point $y$.

\includegraphics[height=6cm]{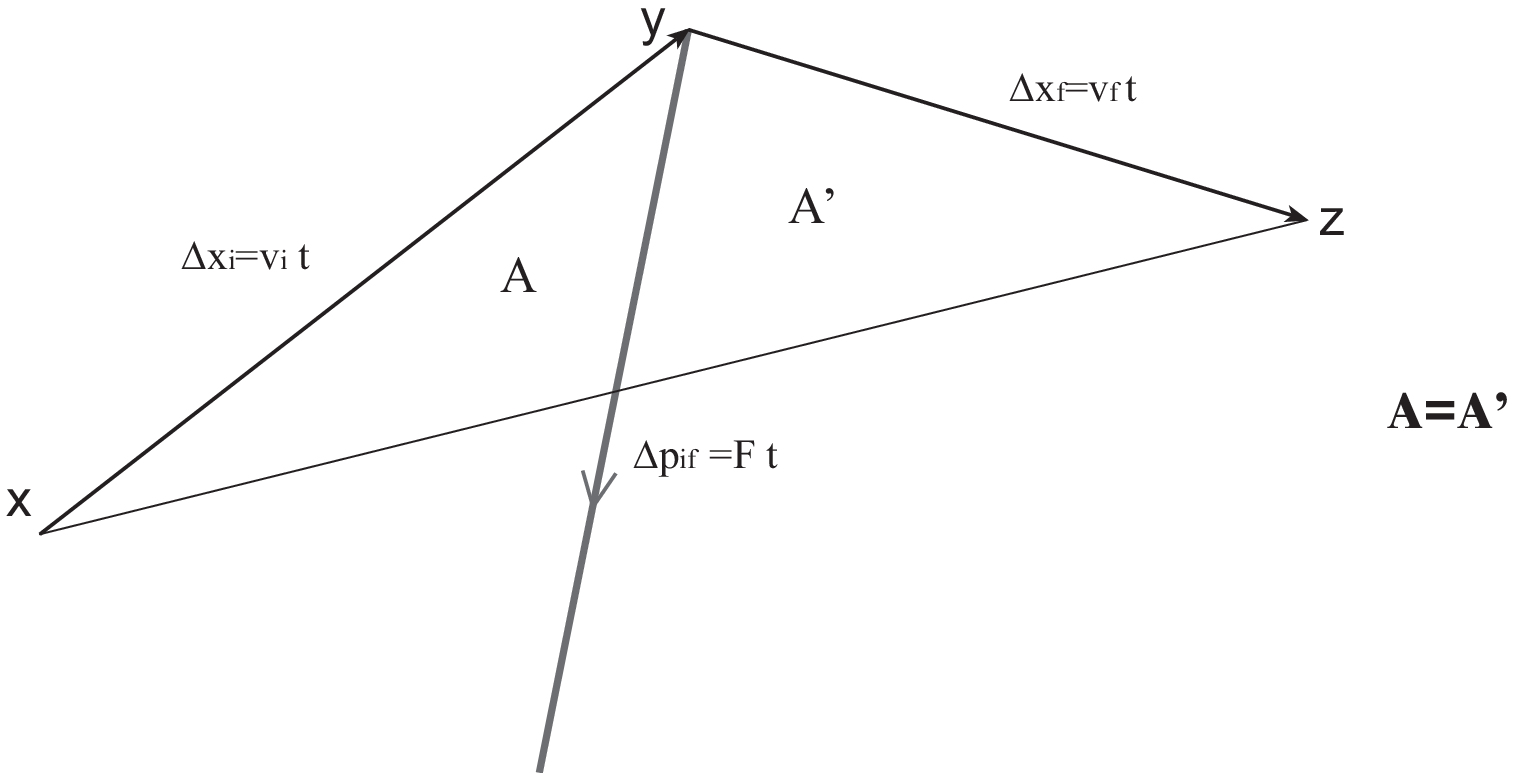}

If we ask for equal-time segments, then areas $A$ and $A'$ are equal too.
In modern language, we are showing that dynamics of a physical system
is given by the rule

  $$  \Delta x_i  \times \Delta p = \Delta x_f \times \Delta p$$

By applying this principle to central forces Newton was able to introduce
time in geometry, mimicking Kepler second law. This is proposition 1 of
Book I in the Principia. Historians tell us that this proposition was 
rebuilt at least three times while doing the built, and it was already
present in the previous paper {\it De Motu}.

It is very troublesome to define evolution, or consistence within a trajectory,
by claiming the equality of two areas, and then asking both areas to go to
zero. Paraphrasing my colleague E. Forgy (from a different context), I believe
the old fathers could be asking themselves: Is mechanics just a series 
of $0=0$ statements? 

Two close remarks:

- It is known that path integral measure is concentrated in continuous
everywhere, differentiable nowhere, trajectories. This shows how troublesome
is to try to approach the classical path. 
By the way, Feynman path integral is about  the limit of an equal-time 
discretization
of trajectory, just as Newton Proposition 1.

- It\^o's stochastic calculus includes a factor $\sqrt t$. It should be 
interesting, from the point of view of didactics, to look for the geometrical
origin of this root. 

\section*{Flash two.}

Consider the natural "elementary school" groupoid over configuration
space:

$$ (x \   y \   u) \circ  (y \  z \  v) = (x \  z \  u+v) $$

Its algebra of functions has the product

$$ (AB)(x,z,t) = \int \int  A(x,y, r) B(y,z,t-r)  dy dr $$

Fourier transforming the $t$ component we see this product is
equivalent to:

$$ (\tilde A \tilde B) (x,  y,  \hat t) = \int \tilde A (x,y,\hat t)
                      \tilde B(y, z, \hat t) dy $$

which in turn, by changing $\epsilon \sim  1/\hat t$, corresponds to
the product defined in the $\epsilon \neq 0$ part of Connes's Tangent
Groupoid. 
Thus the algebra of functions in the latter is a subalgebra
of the one of the "elementary school" groupoid.

The $\epsilon =0$ part of the Tangent Groupoid is defined using the
product   $(x,X) \circ (x,Y) =(x, X+Y)$ of elements of the Tangent Bundle.
  It is well known that continuity of functions in all the groupoid
is just a (de)quantization condition, as $\epsilon >0$ defines a product
of operators in Hilbert Space.

By Fourier transform we can see that Groupoid algebra in Tangent
Space is pointwise product of functions in Cotangent space. And
if we choose a concrete limiting procedure then we can associate a
deformed, star product to the functions of the Cotangent Space. Whose
star-exponential, Fourier transformed again, can be proved to
be Feynman path integral.

\section*{Flash three.}

While the law of areas can be used to prove angular momentum preservation,
its origin is deeper, and simpler, than our modern Noether's theorem. 
It comes from
the combination of the law of addition of vectors and the first law
of Newton, and in this way it simply expresses the equality of the projection 
of velocities in
the plane orthogonal to the variation of momentum. Thus, force -or variation
of momentum-, which describes a plane, is compulsed to be a covector. And
then we can look for the potential functions whose gradient is the force.

The visualization of the plane associated to the force covector let us to
extend the equal-time area law to more arbitrary time steps: now we just
ask for inertia law in [the projection on] this plane. Noncommutativity
is still there, hidden in the definition of variation of momentum just as
something with happens "between" velocity steps (whose size does not 
matter anymore). Lets say, when position changes there is not such thing as
a change of momentum, and when momentum  changes there is not such thing
as a change of position.

It could be interesting to review the introduction of a gauge field. Then
it is known that the new momentum (as defined from the Lagrangian) is not 
anymore the canonical conjugate of position. But the point is that the new
field, at least in the simplest U(1) case, generates a force orthogonal to
the velocity. In some sense it is  exploiting a hole in the "projected inertia
law" of the previous paragraph. One could relate this exploit to the proof
of Maxwell's [homogeneous] laws from that famous -and slippery- report 
of Feynman to Dyson. Remember that noncommutativity is needed for the proof. 

{\footnotesize
\section*{Remark.}

If we were used to the derivation of the Lagrangian (and the Euler-Lagrange
equations) in a deductive way from Newton dynamics, as historically was
done, then neither the groupoid formalism nor Feynman path-integral approach
should surprise us. It is very easy to attach a Lagrangian action to
an element of groupoid and then to formulate the variation mechanism. And
the extremal condition can be shown to be the one of the path integral. More
on this below in the bibliography.

The real quest is the use of these formalisms for fields beyond 0+1.
}

\section*{Flash four.}

During a quantization process it always happens that we choose
an ordering of operators. Above it was the limiting procedure. 
It is equivalent to select one concrete star product. Also, it is
known that it corresponds in the path integral formalism to a choosing
of what discretization method do we apply to the Lagrangian.  

The classical limit of course, does not see the discretization
nor the ordering.

Now, A huge -for the standards of the family- group of theorists like
to study the noncommutative differential calculus coming from the rule
$$f (x) . dx = dx . f(x-\lambda), $$
that approaches the usual calculus when $\lambda \to 0$. 

The previous formula is a bit formal, and it is usually supplemented
by choosing whether the associated difference equations are to be 
considered with $x \in \mathbb{Z}$ or with $x \in \mathbb{R}$. In
any case, there is still the same question that in the physical
methods: While keeping with the noncommutativity rule,
 it could be possible to do various combinations of
forward and backward derivations or to choose different values
of a displacement parameter such that $\lambda f'(x)=
f(\lambda+x+\mu)-f(x+\mu)$ with $\mu\to 0$ when $\lambda\to 0$.


As other possibility,
a doubly twisted differential calculus can be defined via
$$
D_{qr}f(x) ={ f(q x) - f(r x) \over (q-r) x} 
$$
and the standard technologies
$$
D_{qr} x^n = <n> x^{n-1}, \  <n>=
 {q^n-r^n \over q-r} = \sum_{i=0}^{n-1} q^{n-i} r^i , \ \mbox{etc.}
$$

It is interesting because from the physics side
no results are published, as far as I known, beyond $0+1$, so the
q-geometry approach could contain some surprises. 
 It should be good to be aware that quantum mechanics does not need
renormalization, and renormalization is about scale-keeping in a
limit process. 

And even in the trivial case it can be studied an "angle" $q/m$ depending
of the position $x$ but disappearing in the commutative limit.

\section*{Flashback.}
It was not the first time Newton ran into troubles with the infinitesimal
limit. Back in  1666 he noticed how his binomial was the key for an 
algebraic approach to the method of calculation of tangents. The combination
with Barrow theorem, giving the inverse operation, was to be enough to 
dominate analysis during three centuries. Probably Barrow got enthusiastic
by the discovery, because Newton was asked to write a small document, which
was showed to Collins. 
 
It seems that objections were raised about the method, and Newton was
convinced to bury it in the coffin of unspeakable resources.
At that time Barrow influence was high. He had narrowly
 escaped death 
two times while adventuring in foreign lands, trying to rebuild the methods
of Archimedes, and then returning Cambridge to be awarded the (first) Lucasian
chair. And Collins was also remarkable, he was the key editor in the age, even
if his own math was not so impressive. So it is not
surprising that Newton was refrained of speaking calculus, at least
until Leibnitz happened to claim the same results. 

With time -and this is a sad history- Barrow's star should decline,
retired to monastic studies far from math. He evolved to write Collins that
their superiors -in his Order- had completely "forbidden" him to do study 
in any mathematical research, 
and it is rumored he
eventually overdoped opium until death. He was not there, then, to discuss with
Newton about the Principia.

The limit of Newton's Proposition 1 was ultimately justified in empirical
grounds: it fitted the known dynamics of physical bodies. It should take
centuries to verify that the fit was wrong.

\end{document}